\documentclass[12pt]{article}
\usepackage{amsmath}
\usepackage{amssymb}
\usepackage{amsthm}
\usepackage{amscd}
\usepackage{amsfonts}
\usepackage{dsfont}
\usepackage{cases}
\usepackage[applemac]{inputenc}
\usepackage{enumerate}
\usepackage[all]{xy}
\usepackage{esint}
\usepackage{graphicx}
\usepackage{url}
\theoremstyle{plain}

\theoremstyle{definition}

\newcommand{\N}{\mathbb N}

\newcommand{\R}{\mathbb R}

\newcommand{\E}{\mathrm E}
\newcommand{\Var}{\mathrm Var}

\newcommand{\one}{\mathrm\bf 1}

\newcommand{\I}{\mathbb I}

\begin{document}
\title{The BRS-inequality and its Applications\\ A Survey}
\date{}

\author{ F. Thomas Bruss \\Universit\'e Libre de Bruxelles}
\maketitle
\begin{abstract} \noindent This article is a survey of results concerning an inequality, which may be seen as a versatile tool to solve problems in the domain of Applied Probability.
The inequality, which we call BRS-inequality, gives a convenient upper bound for the expected maximum number of non-negative random variables one can sum up without exceeding a given upper bound $s>0.$ One fine property of the BRS-inequality is that it is valid without any hypothesis  about independence of the random variables. Another welcome feature is that, once one sees that one can use it in a given problem, its application is often straightforward or, not very involved.

This survey is focussed, and we hope that it is  pleasant and inspiring to read. The focus is easy to achieve, given that BRS-inequality and its most useful versions can be displayed in  three Theorems, one Corollary, and their proofs. We try to do this in an appealing way. The objective to be inspiring is harder, and the best we can think of is offering a variety of applications. Our examples include comparisons between sums of iid versus non-identically distributed and/or dependent random variables,  problems of condensing point processes, subsequence problems, knapsack problems, online algorithms, tiling policies, Borel-Cantelli type problems, up to applications in the newer theory of resource dependent branching processes. 

Apart from our wish to present the inequality in an organized way, the motivation for this survey is the hope that interested readers may see the potential of the inequality for their own problems. This would contribute to the growing basket of applications, and possibly add new insights. 
\medskip


\medskip
\noindent {\bf AMS 2010 Math. subj. classific.}: 60-01, secondary 60-02.

\medskip\noindent{Short Running title}: BRS-inequality
\end{abstract}

\section{Introduction}

Let $X_1, X_2, \cdots$ be a sequence of positive random variables (possibly dependent) that are jointly continuously distributed and which all have the same absolutely continuous marginal distribution function $F$. For $n\in \{1, 2,\dots\}$ and $s\in \R^+$ let $N(n,s)$ be the maximum \emph{number} of observations from the set $\{X_1,X_2, \cdots, X_n\}$ that one can sum up without exceeding $s.$ More generally, we will study as well the case where the $X_k$ are not identically distributed, and where each $X_k$ has an absolutely continuous marginal distribution function $F_k$. 

The focus of this paper is to find a good upper bound for the expected value $\E(N(n,s))$, and to show its versatility for application.
The inequality which yields such an upper bound is what we call the {\it BRS-inequality.} The reader may find a memory aid  for the meaning of BRS by reading it as 
\medskip

\centerline{\it Budget-Restrained-Sum.}

\medskip\noindent Indeed, one may interpret 
$\{X_1,X_2, \cdots, X_n\}$ as random prices for $n$ items, all having the same utility, offered for sale, and may ask what is the maximum number of them one can expect to be able to buy with a budget limit $s.$

\smallskip
 Clearly one cannot do better than looking at the complete list of prices and first select the smallest price, then add the second smallest price, and so on, as long as the current accumulated bill does not exceed $s.$
Hence we can define $N(s,n)$ equivalently via the increasing order statistics of $\{X_1, X_2, \cdots, X_n\}$, denoted by $X_{1,n} \le X_{2,n} \le \cdots \le  X_{n,n}$, namely
\begin{align}\label{eq:threshold def}
N(n,s)=\begin{cases} 0, {~\rm if}~ X_{1,n}>s,\\
\max\{k \in \N: X_{1,n}+X_{2,n}+ \cdots +X_{k,n} \le s\}, \,{\rm otherwise.}
\end{cases}
\end{align}

\section{The BRS-inequality}
The following Theorem 1, and its generalization Theorem 2, give a simple and useful bound on the expectation of the random variable $N(n,s)$ defined in (1). For Theorem 1 we  add an almost-sure convergence result for the case that the $X_k$ are also independent:
\bigskip

\noindent {\bf Theorem 1} 
~Let all $X_1, X_2, \cdots, X_n$ be identically distributed non-negative random variables with absolutly continuous distribution function $F.$
\smallskip

{\rm (i)}~First, one has
\begin{align}
\E(N(n,s))\le nF(t),
\end{align}
where $t:=t(n,s)$ solves the equation
\begin{align}
~n \int_0^t x dF(x)= s.
\end{align}

{\rm(ii)}~Second, if the random variables are {\it independent} and 
$s:=s_n \to \infty$ in such a way that $\lim_{n\to\infty} n^{-1} s_n$ exists, then 
\begin{align} n^{-1}N(n,s_n)/F(t(n,s_n))~ \to~1~\rm{a.s.}, \text{ as } n \to \infty.
\end{align}

\medskip
\noindent
The equation (3) is visibly central for the BRS-inequality (2), and it is appropriate to call it the  BRS-equation associated with the BRS-inequality, or simply {\it its} BRS-equation. 

\smallskip

The first proof of Theorem 1 with Part (i) and Part (ii) appeared in Bruss and Robertson (1991). We will use Part  (ii) of Theorem 1 later on in several places to strengthen certain results related with Part (i). Our main interest concentrates on Part (i) itself, however. The latter was elegantly generalized in Steele (2016), accommodating non-identically distributed random variables:

\medskip\smallskip

\noindent {\bf Theorem 2} {(Generalized BRS-inequality)}
 \noindent Let $X_1, X_2, \cdots, X_n$ be positive random variables that are jointly continuously distributed and such that each $X_k$ has an absolutely continuous distribution function $F_k.$ If $N(n,s)$ is defined as in (1), then
\begin{align}\E(N(n,s)) \le \sum_{k=1}^n F_k(t),
\end{align}where
$t:=t(n,s)$ is the unique solution of the equation
\begin{align}\label{eq:new def of t(n,s)}
\sum_{k=1}^n \int_0^t x dF_k(x)= s.
\end{align}

\medskip\noindent Clearly, if all random variables $X_i$, $i \in\{1,2,\cdots, n\}$, have the same marginal
distribution $F$, then (6) recaptures (3), and (5) recaptures (2). 

\medskip\noindent
The proofs of Theorem 1 (i) and (ii), and Theorem 2, are thus in two different journals, and our first goal is to offer here a single combined proof combined with  relevant remarks in the view of possible applications. The second objective is to offer with a new theorem (Theorem 3) a refined BRS-inequality and also to briefly discuss other possible directions of refinements. Our third goal is to give several examples from different domains of Applied Probability, showing the versatility of the inequality for applications.

\bigskip
\noindent {\bf Proof}~ We start with the proof of Part (i) of Theorem 1, and we do this in a way that leads almost directly to the proof of the more general Theorem 2.

\smallskip\noindent By the continuity of the joint distribution of $(X_i: 1\le i\le n)$,
there is clearly a unique set of selected indices $A\subseteq \{1,2,\cdots, n\}$ of the variables that attains the maximum
in the definition of $N(n,s)$ in (1).
We denote this subset by $A(n,s).$ Note that $\#A(n,s)=|A(n,s)|=N(n,s).$ 

\noindent Further let \begin{equation}\label{eq: B def}
B(n,s)= \{i \in \{1,2,\cdots, n\}: X_i \leq t(n,s) \},
\end{equation}
where $t(n,s)$ is the threshold value determined by the implicit relation (3).

\smallskip
The idea  is now to compare the sets $A(n,s)$ and $B(n,s)$,
together with their associated sums,
\begin{equation}\label{eq:two sums}
S_{A(n,s)}=\sum_{i \in A(n,s)} X_i \quad \text{and} \quad S_{B(n,s)}=\sum_{i \in B(n,s)} X_i.
\end{equation}
From these definitions, one immediately gets by definition  two useful relations, namely
\begin{equation}\label{eq:two relations}
S_{A(n,s)}\leq s \quad \text{and} \quad \E[S_{B(n,s)}]=n \int_0^{t(n,s)} x \, dF(x) = s.
\end{equation}

Now, by its definition, $S_{A(n,s)}$  is a partial sum of order statistics. The summands of
$S_{B(n,s)}$ consist precisely of the values $X_{n,i}$ with $X_{n,i} \leq t(n,s).$ We can think of the latter as also being listed in increasing order, of course, so that
$S_{B(n,s)}$ is also equal to a partial sum of order statistics of  $\{X_1,X_2, \ldots, X_n\}$.
These observations will help us with  estimations of the relative sizes of the two sums $S_{A(n,s)}$  and $S_{B(n,s)}$.

\smallskip
First, we note that one must have either $A \subset B$ or $B \subseteq A$.
More specifically, if
$S_{B(n,s)}\leq S_{A(n,s)}$ then we have $B(n,s) \subseteq A(n,s)$. Moreover, the summands $X_i$ in the difference set, that is,  with $i \in A(n,s) \setminus B(n,s),$ these summands are all bounded
below by $t(n,s).$ Hence we have the bound
$$
S_{B(n,s)} + t(n,s)\left(|A(n,s)|-|B(n,s)|\right) \leq S_{A(n,s)}, \quad \text{if } S_{B(n,s)}\leq S_{A(n,s)}.
$$
Similarly, if
$S_{A(n,s)} < S_{B(n,s)}$ then $A(n,s) \subset B(n,s)$ and the summands $X_i$ with $i \in B(n,s) \setminus A(n,s)$ are all bounded
above by $t(n,s)$; so, in this case,$$
S_{B(n,s)}  \leq S_{A(n,s)} + t(n,s)\,\big(|B(n,s)|-|A(n,s)|\big), \quad \text{if } S_{A(n,s)}< S_{B(n,s)}.
$$
Taken together, the last two relations tell us that whatever the relative sizes of $S_{A(n,s)}$ and $S_{B(n,s)}$ may be,
one always has the {\it key inequality}
\begin{equation}
t(n,s)\left(|A(n,s)|-|B(n,s)|\right) \leq S_{A(n,s)}- S_{B(n,s)}.
\end{equation}
Here $t(n,s)>0$ is a constant. We recall that $|A(n,s)|=N(n,s)$, and we see from \eqref{eq:two relations} that the right-hand side has non-positive expectation. Hence,
taking the expectations in the key inequality (10)  gives us
$$
\E\left( N(n,s)\right) \leq \E \left(|B(n,s)|\right) =
\E\left(\sum_{i=1}^n \I \left\{X_i \leq t(n,s)\right\}\right) = n\,F(t(n,s)),
$$
and the proof of the bound (2) is complete.

\medskip
Now, to prove its more general form as given in Theorem 2,  we again define the set $B(n,s)$ by (7).
Additional care is here needed with the definition of  $A(n,s)$, and, as we shall see in the subsequent Remark 1, this care  provides at the same time a non-negligible benefit.

\smallskip

Specifically,
we first define a total order on the set $\{X_i: i \in \{1, 2, \cdots, n\}\}$ by writing $X_i \prec  X_j$ if either condition (A) or condition (B) holds, where

\bigskip
(A) ~ $X_i < X_j$, 

\smallskip
(B) ~ $X_i=X_j$ and $i < j$.

\bigskip\noindent
Using this order, there is now a unique permutation $\pi:\{1, 2, \cdots, n\} \rightarrow \{1, 2, \cdots, n\} $ such that
$$
X_{\pi(1)} \prec X_{\pi(2)} \prec \cdots \prec X_{\pi(n)},
$$
and we can then take $A(n,s)$
to be largest set $A \subset \{1, 2, \cdots, n\}$ of the form

\begin{equation}\label{eq:second def of A(n,s)}
A= \{ \pi(i): X_{\pi(1)} + X_{\pi(2)} + \cdots + X_{\pi(k)} \leq s \}.
\end{equation}

\smallskip
We can now proceed with the proof of key inequality (10)
in the same way. We use the identical definitions \eqref{eq:two sums} for the sums $S_{A(n,s)}$ and $S_{B(n,s)}$, because these definitions are not affected by the fact that the variables may now stem from different distributions. Hence, 
by the new definition \eqref{eq:new def of t(n,s)} of $t(n,s)$, we now have
\begin{equation*}\label{eq:E of S sub B}
\E\left(S_{B(n,s)}\right)=\sum_{i=1}^n \int_0^{t(n,s)} x \, dF_i(x) = s.
\end{equation*}
Since we still have $S_{A(n,s)}\leq s$, the expectation on the right side of inequality (10) is again non-positive. 

But then, since all the following arguments in the proof of Part (i) of Theorem 1 are valid independently of the number of different distributions functions $F_k$ which are involved, 
the proof of Theorem 2 is also complete.\qed

\bigskip
As the organization of the argument makes explicit, none of the required calculations
require more on the joint distribution of the variables  $X_i: i \in\{1, 2, \cdots , n\}$ than the joint continuity.
The argument just uses point-wise bounds and the linearity of expectation.

\subsubsection{Terminology}
In the following we will call invariably (2), respectively (5), the BRS-inequality, and (3), respectively (6), the corresponding BRS-equation. It will always be evident from the context of the problem
when (6) and (5) are meant. It is also convenient to always speak of {\it budget} if we refer in the text to the values $s$ or $s_n$ in  these equations or inequalities. 
Moreover, as announced already,  we will also present two refinements of Theorem 2. These refinements are in the same spirit as Theorem 2, and, in particular, they do not require to revise  the definitions of (5) and (6) in any way. Therefore we may and will call them, without a risk of any confusion, also BRS-inequalities.

\subsection{Independence and almost-sure convergence}

\noindent The asymptotic result (4) is proved in
section 2 of the paper Bruss and Robertson (1991).  Note that the general setting of Theorem 2 does not allow to include (4) directly, where independence of the variables is needed. 

The convergence in (4) in Part (ii) has its interest on its own as we shall see later on in several examples.  Therefore we recall here the proof, but in a slightly different way, namely in terms which make the result very intuitive:
\medskip

\noindent {\bf Proof} (Theorem 1, Part (ii)).

\smallskip
\noindent First we note from the definition of $N(n,s)$ in (1) that, although the random time $N(n,s)$ is not a stopping time with respect to the usual filtration generated by the random variables $X_j$, the random time $N(n,s)+1$  {\it is} a stopping time on the filtration $(\tilde {\cal F}_k)$, where $$\tilde {\cal F}_k=\sigma(X_{1,n}, X_{2,n}, \cdots X_{k,n}), k=1, 2, \cdots n.$$ If we sum up all those observations which are smaller than  $t_n$, say, then Wald's equality can be applied for the expected stopping time since, in Theorem 1,  all variables follow the same distribution and thus have the same mean. Moreover, by independence, the unordered variables not exceeding $t(n,s)$ are continuous i.i.d random variables on the interval $[0,t(n,s)].$  
Letting $n\to \infty,$  the Strong Law of Large Numbers implies then that the average of these terms satisfies $$ \frac{1}{n}\,\sum_{k=1}^n X_k\,\one\{X_k \le t(n,s)\} \to \int_0^{t(n,s)} xdF(x)~{\rm a.s.~~as}~n \to \infty,$$ that is, the sum of these selected variables converges almost surely to the truncated expectation 
figuring in (3).  

Finally, if $s:=s_n$ grows with $n,$ and if we let now the threshold become $t:=t(n,s_n)$ then the same argument remains valid provided that $s_n/n$ tends to a limit. Then (3) yields (4),  and  the proof is thus complete.

\qed

\bigskip\noindent{\bf Remark 1~}
Bruss and Robertson (1991) used a setting of i.i.d. random variables because they were motivated by problems in which the result (4) played the main role.  They saw that  independence was not needed for proving (2) with (3) in Theorem 1 but they did not draw attention to this fact.  As Steele (2016) pointed out, the inequalities gain much interest by knowing that no independence is needed. This is also the main reason for introducing the total order induced by (A) and (B). Indeed, note that this allows us to have variables $X_j$ and $X_k$ for $k\not= j$ coincide with positive probability.

\subsection {Refining the BRS-inequality}

\smallskip
It is straigthforward to refine the BRS-inequality by fully exploiting the key inequality (10).
This will be quickly done, and the reader may ask why we did not do this right away before. We will shall shortly discuss the reason  after the Theorem.

\bigskip\noindent{\bf Theorem 3~} Under the conditions of Theorem 2 we have
\begin{align*} \E(N(n,s)) \le \sum_{k=1}^n F_k(t(n,s))-\frac{s-\E(S_{A(n,s)})}{t(n,s)},
\end{align*}
where $t(n,s)$ satisfies (6), and where $S_{A(s,n)}$ is the sum defined in (8).

\bigskip\noindent{\bf Proof}
~ The key inequality (10) can be written in the form
\begin{align}
| B(n,s)|\,-\, |A(n,s)|\,\ge \frac{S_{B(n,s)}-S_{A(n,s)}}{t(n,s)},\end{align} where $t:=t(n,s)$ solves equation (6).  Recall again that $S_{A(n,s)}$ defined in (8) is the sum of the $N(n,s)$ smallest order statistics, and thus $|A(n,s)|=N(n,s).$ Further, since $S_{B(n,s)}$ is the sum of all those variables not exceeding $t(n,s),$ we have
\begin{align} |\,B(n,s)\,| = \sum_{k=1}^n \,\one\{X_k\le t(n,s)\} \rm{~and}~~ \E(S_{B(n,s)})=s, ~\end{align} where the second equality was shown in  (9) for identically distributed random variables, and for different distribution functions correspondingly after (11).

Using these equalities we obtain from (12) by taking expectations on both sides \begin{align}\left( \sum_{k=1}^n F_k(t(n,s))\right)-\E(N(n,s))\ge \frac{s-\E(S_{A(n,s)})}{t(n,s)}. \end{align}
Substracting in the last inequality the sum term on both sides, and then changing signs, yields the statement Theorem 3. 
\qed

\bigskip\noindent 
{\bf Remark 2}
~We are of course aware that Theorem 1 and Theorem 2 are corollaries of Theorem 3, and that, following our tradition in Mathematics,
Theorem 3 should be entitled to be stated and proved first. Note however that, if we had done so,  it would have been necessary to explain the meaning of $A(n,s)$ and $S_{A(n,s)}$  already in the hypotheses of the  theorem. This would look technical and probably distract from the main interest of the inequality. Therefore the reader may agree that, viewing readability, the way we ordered the theorems is nicer.  

\subsection{\bf Interest of refinements}
\medskip
Theorem 3 is of interest as soon as we have additional information. 

\smallskip\noindent 
If we call the random quantity $s-S_{A(n,s)}$ the {\it residual} (of available budget), and if we have sufficient information to give easy bounds for its expectation $\E(s-S_{A(n,s)})$, then Theorem 3 can of course give us a more precise bound. This is (independently) true if $s$ and $n$ are such that $P(N(n,s)<n)$ is very close to $1$ so that the $S_{A(n,s)}+X_{N(n,s)+1, n}-s$ represents the (very probable) {\it over-shoot}. In specific examples, such as for instance when all $F_k$'s have similar bounded supports, Theorem 3 helps us also to understand why the bound we obtain  Theorems 1 and 2 are  often quite tight.
\medskip

The latter type of arguments allows us also to defend the attitude that any strengthening the BRS-inequality is always of interest, at least  in some special cases. Using standard arguments about conditional expectations, we obtain a version of the BRS-inequality which, with sufficient information about (nice) properties of the $F_k$ and their supports,
becomes  "tight".  

\smallskip
The following inequality is such a version of the BRS-inequality. It is a corollary of Theorem 3.

\bigskip\noindent
{\bf Corollary} Let $N(n,s)$ be defined as in (1), and $t:=t(n,s)$ as in (6). Also let the set $A(n,s)$ be defined by (8) and recall that $\E(S_{A(n,s)})\le s.$ Denoting $p_n=P(N(n,s)<n)$ we have

\begin{align*} E(N(n,s))
\le p_n\left( \sum _{k=1}^n F_k(t)-\left(\frac{s-\E(S_{A(n,s)})}{t}\right)\right) + n(1-p_n).\end{align*}

\smallskip \noindent {\bf Proof}~~By using the evident equality $E(N(n,s)|N((n,s)\ge n)=n$ and conditioning on the event $\{N(n,s)<n\}$ and its complement, we obtain 
$$\E(N(n,s))= n (1-p_n)+p_n\, \E\left(N(n,s)\Big| N(n,s)<n\right).$$ The Corollary follows now in a straightforward manner from the trivial inequality $$
E(N(n,s)\Big|N((n,s)<n)\le \E(N(n,s)$$
by plugging in the upper bound for $\E(N(n,s)$ of Theorem 3, and then isolating $\E(N(n,s)$ again on the left-hand side.\qed

\bigskip\subsection{Interest of further refinements}

It is now no surprise that, in some special cases, further refinements should be possible. We want to exemplify another one of such cases, but  we confine ourselves to describe things verbally, because it suffices to present ideas which come along quite naturally. There are simply too many special cases which come to our mind.

\smallskip
Consider for instance the case of identically distributed random variables  $X_j\sim X$ (as in Theorem 1). Suppose moreover that $s:= s_n$ is small compared to $n\E(X).$ Then, for large $n$, the probability of not exceeding the budget $s_n$ becomes small, and so $p_n=P(N(n,s)<n) $ becomes large as $n$ grows.  Consequently, it will become very probable that the budget will be exceeded (by definition  upon observation of $X_{N(n,s)+1,n} \in\{X_1, X_2, \cdots ,X_n\},$ and we know  $$P(X_{N(n,s)+1,n} \in\{X_1, X_2, \cdots ,X_n\}) \to 1 ~~{\rm as }~n\to \infty.$$
This guarantees that the sum $S_{A(n,s_n)}$ behaves essentially like $s_n$ as $n$ becomes large.
Now recall that, conditioned on $N(n,s)<1,$ the 
random time $N(n,s)+1$ is a stopping time with respect to the filtration $$\tilde {\cal F}_k=\sigma(X_{1,n}, X_{2,n}, \cdots X_{k,n}), k=1, 2, \cdots, n.$$  Recall also that  all $X_j \sim X$ have the same mean and that the generalized version of Wald's equation works with weaker conditions than independence. Thus, even though we started with a special case and, compared to independence, the incomparable (arguably partially weaker) condition that $s_n$ is small compared to $n\E(X),$ this allows then a similar approach than shown in Subsection 2.1.

\bigskip
Having said this, the author thinks that in this survey it would be hardly rewarding to go further
with discussions of more specific cases. We probably should give priority to confining ourselves to results 
which are, arguably, of sufficient general interest.
 If a researcher would like to apply a result for a given problem, needing slightly more precise estimates, and he or she knows of the existence possible refinements, then she or he, knowing their own problem best, would also know best in which corner of the set of options to look for. 
 As so often in Mathematics, it is not always the very strongest form of a result which draws attention but rather the most convenient form of a sufficiently strong result, which does so. 

Seen under this angle of view, our analysis of relative importance is as follows.
\subsection{Evaluation of importance}
Theorem 1 and its generalization Theorem 2 give a convenient upper bound of $\E(N(n,s)).$ They involve those functions and parameters which intervene naturally, namely the number of random variables $n$, the available space or budget $s,$ and the distributions of the variables. They neither speak about specific properties of the latter nor would they refer indirectly to such properties.   It is this "modesty" in the hypotheses that makes, as we think, Theorems 1 and 2 appealing for a wide range of applications. 

Theorems 1 and 2 are in that sense likely to be the most interesting results.  The author thinks that Theorem 3 is of sufficiently general interest, but that its corollary is already somehow close to the borderline.  

As we will see later, Theorem
1 has already a considerable number  of citations to its credit. Thus so far it stands out. The likely reasons are easy to understand. On the one hand, the univariate case is typically more frequent in studied problems than a case involving several different distributions. On the other hand, and in particular, Theorem  1 (Bruss and Robertson (1991)) has been around much longer than the (partially) generalized version Theorem 2 (Steele (2016)).

\smallskip 

The author thinks that Theorem 2, i.e. Steele's version of 2016, is the most promising from the point of view general interest for future applications. Theorem 3 is also very welcome on this preference list of the BRS-inequality, since there are many applications  where additional information comes in naturally, and then Theorem 3 sharpens results for both Theorem 1 and Theorem 2.

\section{Applications}
The goal of this section, as announced before, is to convince that the BRS-inequality is an interesting result, and both flexible and versatile for applications.

\smallskip \noindent We begin our examples with identically distributed uniform random variables.

\subsection {Identically distributed $U[0,1]$ variables}

\noindent Let $X_1, X_2, \cdots, X_n$ be uniform random variables on $[0,1],$ that is $F(x)=x$ on $[0,1].$ To compute an upper bound
for the maximum number of $X_k$'s we can sum up without exceeding $s$ with $0<s<n/2$ we solve the BRS-equation (3), that is
\begin{align} n \int_0^{t(n,s)} x dF(x)\,=\,n \int_0^{t(n,s)} x dx\,=s\, ~\iff~ t(n,s)=\sqrt{2s/n}.\end{align} Hence, from the BRS-inequality \begin{align}\E(N(n,s))\le n\,F\big(\sqrt{2s/n}\big)=\sqrt{2sn}.\end{align} This bound becomes trivial for $s\ge n/2=n\E(X),$ namely $\E(N(n,s))\le n.$ We can say slightly more, i.e. $\E(N(n,n/2)<n,$ because $n$ items need in expectation the budget $n/2.$ However, they need  a larger budget with a strictly positive probability.

If the $X_k$'s are moreover independent, then  this bound is essentially tight according to (4).

\subsection{ Random variables with dependence}

\smallskip\noindent Recall that  no independence assumption about the $X_k$'s is needed. So let us compare the extreme case of total independence with the case of full dependence such as $X_1=X_2=\cdots =X_n$, say. Here, by the way, we see that the total order defined in the proof before (11) through (A) and (B) is essential. In the case of uniform random variables on $[0,1]$ we know from Subsection 3.1 that $\E(N(n,s))\le n F(t(n,s))=\sqrt{2sn}$.  Moreover in the  case of i.i.d. $U[0,1]$ random variables we have from (4) also the asymptotic relationship $N(n,s_n)\sim \sqrt{2s_n}$ so that the upper bound computed in (15) and (16) is asymptotically tight. 

In the case $X_1=X_2=\cdots =X_n$ however, there is trivially only one order statistic, namely $X_1, $ and our bound should be worse. Indeed, it is easy to compute it precisely, namely $$\E(N(n,s))=\sum_{k=1}^nP(N(n,s)\ge k)=\sum_{k=1}^nP(X_1\le s/k)$$$$
=\sum_{k=1}^n\left(\one\{k<s\}+s \one\{k\ge s\}/k\right) \sim [s]+s \log(n),$$ where [s] denotes the largest integer smaller than $s.$ This upper bound grows more slowly in $n$ than $\sqrt{2sn},$ and it is clear why: independent draws bear naturally the potential of producing
small order statistics, whereas these independent draws are now wiped out by the total dependence. 
\subsection{Different distribution functions $F_k.$}

\smallskip \noindent For each $k \in \{1,2, \cdots,n\}$ let $X_k$ to be uniform on the real interval $[0,k].$
For convenience we take $0 < s \leq 1$ and $n \geq 4.$ The BRS-equation tells us
\begin{equation}\label{eq:determineS}
s=\frac{1}{2} \sum_{i=1}^n \frac{1}{i} \tau^2(n,s)=\frac{1}{2} t^2(n,s) H_n \quad \text{or} \quad  t(n,s)=\sqrt{2s/H_n},
\end{equation}
where $H_n$ denotes the $n$'th harmonic number.
In particular, for $s=1$  \begin{equation*}
\E(N(n,1))\leq \sum_{i=1}^n F_i(t(n,1))=\sum_{i=1}^n \frac{1}{i}({2}/{H_n})^{1/2}
=(2 H_n)^{1/2},
\end{equation*}
where we use $n\geq 4$ to assure that $H_n> 2.$ 
This bound grows also much more slowly than $\sqrt{2n}$ computed in 2.1. The reason is now that small order statistics become rarer with the increasing supports $[0,k]$ of the uniform $F_k$'s. 

\smallskip
If, in contrast, we have decreasing supports as $k$ grows, the upper bound for $\E(N(n,s))$ should go up. We would expect the upper bound to become trivial or close to trivial if the supports of the $F_k$ are decreasing more or less in the same order (or quicker) as/than the expected order statistics of $n$ variables would do on the support of $F_1.$ To exemplify this, let us choose the {\it exact} corresponding order, i.e. $$F_k(x)= \I\{0\le x\le 1/k\}kx+\I\{x>1/k\}, k=1, 2, \dots, n,$$ and again $s=1.$ It is easy to check that the BRS-equation becomes
$$\sum_{k=1}^n k \min\{1/k,t\}=1=\sum_{k=1}^n~kt=t\,\frac{n(n+1)}{2},$$ where we used that, in the first equality, we must have $t:=t(n,1)\le 1/n,$ because the sum cannot equal 1 otherwise. Hence the solution is $t=2/(n(n+1)).$ Plugging $t$ into the BRS-inequality  (5) yields  $\E(N(n,1))\le \sqrt{n(n+1)}\sim n,$ showing that our bound is here almost perfect.\bigskip
\subsection{Point processes and selection bias.}

\smallskip \noindent
Let $T_0=0$ and, recursively, $T_k=T_{k-1}+X_k, k=1, 2, \cdots$ where the $X_k$'s are positive random variables with absolutly continuous laws $F_k.$ The $X_k$'s are now seen as the inter-arrival times of a point process $(T_k)_{k=1, 2, \cdots}$ on the positive half-line.

Suppose now that we fix $n$ and wait for the $n$th arrival,  occurring at time $T_n=T$, say.  Given $s$ with $0<s<T$ we ask what is the maximum density of arrival points we can expect on a set ${\cal S} \in U[T]$ where $U[T]$
is the set of all subsets of $\{[T_{k-1},T_k[: 1 \le k \le n\}$?

The maximum obtainable density in this definition is simply $N(n,s)/s$ with $N(n,s)$ defined in (1) based on the inter-arrival times $X_1, X_2, \cdots, X_n.$ Thus from the BRS-inequality we know  that the  maximum expected density is  bounded by \begin{equation}d_{\max}(s):=\frac{\E(N(n,s))}{s}\le \frac{1}{s} \sum_{k=1}^n F_k(t(n,s)),\end{equation} where  $t(n,s)$ solves the BRS-equation (3) if the $X_k$ are identically distributed, and otherwise solves (6). We also note that in the case of independence (where we have  the stronger result (4)), the problem of finding the set
with {\it minimum} expected density is then {\it asymptotically} dual, because it suffices to maximize the number of points on the complementary set with Lebesgue measure $T_n-s.$ (see Figure 1.)
\begin{figure}[ht]
\centering\includegraphics[width=0.90\textwidth, angle=0]{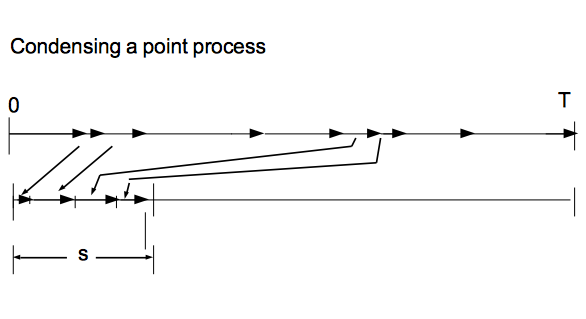}
\centerline{ Figure 1}
\begin{quote} Figure 1 shows for a given "budget time" $s$ the selection of the $4$ smallest inter-arrival times (after having seen them all) until arrival time $T=T_9$.\end{quote}
\end{figure}

\subsection{ Poisson processes}
 As a specific  example, consider a homogeneous Poisson process of rate $1$. Hence here all the $F_k$ are the same.
 
 Suppose that a dishonest statistician would like
to make  an observer believe that the rate is higher than $1$ by claiming that the missing data are due to the fact that the counting process could not be recorded during certain sub-periods. 
How far can this type of dishonesty bias the perceived rate?

\begin{figure}[ht]
\centering\includegraphics[width=1\textwidth, angle=0]{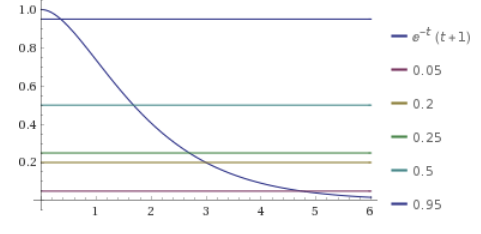}
\centerline{Figure 2}
\begin{quote}Figure 2 displays for a Poisson process of rate $1$ the l.h.s. of (16) for $0\leq t\leq T=n$, and the r.h.s. of (16) for a selection of values of $s/n$ ranging between 5\% and 95\% (horizontal lines $\equiv$ fraction of the observation time.) \end{quote}
\end{figure}

\smallskip 
 
It is well known that the inter-arrival times of the Poisson process (rate $1$) are i.i.d. exponential random variables with distribution function $F(x)=1-e^{-x}.$ No conditional bias intervenes if we look at the random horizon $T=T_n.$ The BRS-equation (3) yields then
$$\int_0^t x dF(x)=\int_0^t x e^{-x}dx= \frac{s}{n},$$ or equivalently, the implicit equation\begin{align} e^{-t}(t+1)=1 -\frac{s}{n}.\end{align}
Thus (15) yields with $t:=t(n,s)$\begin{equation}d_{\rm \max}(s)\le \frac{1}{s} n F(t(n,s))=\frac{1}{s} n \left (1-e^{-t(n,s)}\right), \end{equation} where $t:=t(n,s)$ solves (19).
Here we note that in the Poisson process case the random variables $X_k$ are i.i.d., and so (4) implies that the obtained bounds obtained from the inequality (20)are essentially tight. The $50\%$ line for instance (see Figure 2) cuts the curve in $t\approx1.67835,$ and (19) and (20) yield a maximum density  $\approx 1.62664,$ that is about 63 percent more than the true rate. Clearly, the less (selected) time we observe, the more the observed rate can grow. With only $5\%$ observation time the dishonest statistician can  cheat (a naive client) by a factor almost $6.$

\subsection{  "Awkward" point processes}

\smallskip
~~~~~~Sometimes one comes across random structures which  are complicated, but have components which are nevertheless well understood, and/or tractable.
Staying in the class of problems we have considered in Subsections 3.4.and 3.5, the following "awkward" counting process (completely invented) exemplifies such an instance.
 
\smallskip
Let $(X_k)$ be a random sequence of inter-arrival times of a point process. $X_1$ is uniform on
$[0,1]$, say, $X_2=1-X_1$, $X_3=1-X_2=X_1$, ... and this alternating pattern goes on until some random time $C_1.$  $X_{C_1+1}$ is then again $U[0,1]-$distributed independently, and the following variables play the same alternating game as before until some random time $C_2.$ Then  again a new U[0,1] random variable
$X_{C_2+1}$ is drawn ...  and so on. Suppose we are stopped after a large number of  $n$ observations  (in total).
Can we say anything nontrivial about the expected number of variables among $X_1, X_2, ...., X_n$ we can sum up without exceeding a given $s$?

At first glance it may seem that we must know something about the mechanism producing the random times $C_1$, $C_2, \cdots.$ But then we simply concentrate on the essential information which is given: whatever the law of $X_1, X_2, ...$, whatever their dependence structure,  all the {\it marginals} of the $X_k$ are the same, namely
$F(t)=t$ on $[0,1].$ Hence, as seen in Subsection 3.1, we always have the convenient bound \begin{align}\E(N(n,s))\le \sqrt{2 s n}.\end{align}
If we happen moreover to know one of the expected change times, $\E(C_n)$ say,  then we have also a nice upper bound for $\E(N(C_n,s))$, namely by the tower property of conditional expectations, 
\begin{align*}E(N(C_n,s))=E\big[E(N(C_n,s)\big|C_n)\big]\le\E\left(\sqrt{2s\, C_n}\right)\le\sqrt{2s\, \E(C_n)},\end{align*}
where the first inequality follows from (21), and the second one from the concavity of the square-root function and Jensen's inequality. 
 
 \smallskip
With more information one can sometimes say more. To stay with our example, suppose all block lengths  $C_{i+1}-C_i$'s of strictly dependent random 
variables are geometric random variables with parameter $p,$ say. Instead of looking at the horizon of $n$ variables  we now look at $n$ blocks, that is $C_n$ variables with budget 
$s:=s_{C_n}.$ The expected length of blocks consisting of strongly dependent variables is thus $1/p.$ 

On the one hand, the larger $p$ gets, the higher becomes the fraction of independent $U[0,1]$ random variables so that upper
bound for $C_n$ variables $$\sqrt{2 s_{C_n} C_n}\sim \sqrt{2 s_{C_n}n/p}$$ must become better and better. On the other hand, if $p$ becomes small, the expected block lengths become large so that most blocks contain with each $X$ also the value $1-X.$ This opens then a new option, namely we  can recompute a new threshold by taking the variables two by two.

\bigskip
\noindent{\bf Heuristic arguments}: All what we say here is that in point processes consisting of both i.i.d. inter-arrival times strongly dependent ones, it can be worth checking whether the fraction of i.i.d. random variables is non-negligible, knowing that the BRS-bound is essentially tight for these. For instance, if they are $U([0,1]$ random variables, $N(n,s)$ is in the order of square-root of their number. The contribution of the others will be somewhere between a logarithmic order (see Subsection 3.2) and an order of square-root, and may contribute to $\E(N(n,s))$ not more than  something asymptotically negligible.

\subsection{Sequential Subsequence Selection Problems and Knapsack Problems}

As has been noticed by several authors before, as e.g. Gnedin (1999) and Arlotto et al. (2015a , 2015b), the BRS-inequality leads to \emph{a priori} upper bounds for two well-studied problems in combinatorial optimization.
In particular, in the classical case of independent uniformly distributed random variables,
the BRS-inequality (Theorem 1) gives bounds that are essentially sharp
for both the sequential knapsack problem and the sequential (increasing) subsequence selection problem (Steele (2016)).

\subsubsection {Sequential knapsack problem}
In this problem, one observes a sequence of $n$ independent non-negative random variables $X_1, X_2, \ldots , X_n$
with a fixed, known distribution $F$. One is also given a real value $x \in [0,\infty)$ that one regards as the capacity of a knapsack into which
selected items are placed.  Observations are sequential, and when $X_i$ is first observed, one either selects $X_i$ for inclusion in the knapsack or
else $X_i$  is rejected from any future consideration. The goal is to maximize the expected
\emph{number}  of items, $\tilde N(n,x)$, say of items that can be sequentially packed without recall into the knapsack with initial capacity $x$. Since the BRS-inequality tells us
how well we could do if we knew in advance all of the values $\{X_i: 1 \le i \le n\}$, it is
evident that no strategy for making sequential choices can ever lead in expectation to more affirmative choices than $\E(N(n,x)),$ and thus from (2) we have the upper bound
\begin{align}
\E(\tilde N(n,x)) \le \E(N(n,x))\le n F(t(n,x)),
\end{align}
where $t:=t(n,x)$ solves the associated BRS-equation (3).

Now, packing an item into the knapsack will not change the distribution of the items to come. The only effect of packing an item of size $y$ say is reducing the current remaining capacity by $y$. Hence the problem is a Markov decision problem so that the optimal sequential selection strategy given
by a unique non-randomized Markovian decision rule. Beginning with $n$ values to be observed and
with an initial knapsack capacity of $x$, the expected number of selections that one makes is denoted by $v_n(x)$. It is easy to see that the value function for this Markov decision problem can be calculated by the recursion relation
\begin{equation}\label{eq:Bellman Knapsack}
v_n(x)=(1-F(x)) v_{n-1}(x) + \int_0^x \max \{ v_{n-1}(x), 1+ v_{n-1}(x-y) \}\, dF(y).
\end{equation}
Specifically, one begins with the obvious relation $v_0(x) \equiv 0$, and one computes $v_n(x)$ by iteration.
This is the Bellman equation (optimality equation)
for the sequential knapsack problem.

\subsubsection {Monotone subsequence problem}
Now suppose that we observe sequentially $n$ independent
random variables $X_1, X_2, \\\cdots, X_n$ with the common continuous distribution $F.$ Our goal is to make monotonically (decreasing, say) choices
$$
X_{i_1} > X_{i_2} > \cdots > X_{i_k},
$$
trying to maximize the expected
number of choices. 
Here, prior to making any selection, we take the state variable $x$ to be the supremum of the support of $F$, which may be infinity. After
we have made at least one selection, we take the state variable $x$ to be the value of the {\it last} selection that was made.
Now we write $\widetilde{v}_n(x)$ for the expected number of selections made under the optimal policy
when the state variable is $x$ and where there are $n$ observations to come.
In this case the Bellman equation given by
Samuels and Steele (1981) can be written as
\begin{equation} 
\widetilde{v}_n(x)=(1-F(x)) \widetilde{v}_{n-1}(x) + \int_0^x \max \{ \widetilde{v}_{n-1}(x), 1+ \widetilde{v}_{n-1}(y) \}\, dF(y),
\end{equation}
where again one has the obvious relation $\widetilde{v}_n(x) \equiv 0$ for the initial value.
In (20)
the decision to select $X_1=y$ would move the state variable to $y$, so here we have
$1+\widetilde{v}_{n-1}(y)$ where earlier we had the term  $1+ v_{n-1}(x-y)$ in the knapsack Bellman equation.
In the knapsack problem the state variable moves from $x$ to $x-y$ when $X_1=y$ is selected.

In general, the solutions of (23) and
(24) are distinct, but Coffman et al. (1987) made the essential observation that $v_n(x)$
and $\widetilde{v}_n(x)$ are equal when the observations are uniformly distributed. This allows to create an interesting detailed linkage between these two problems: 
The first step is to note that the equality of the value functions
permits one to construct optimal selection rules that can be applied simultaneously to the same sequence of
observations. The selections that are made will be different in the two problems, but we see a useful
distributional relationship.
The essential observation is that the second term of the
Bellman equation leads one almost immediately to the construction of an optimal
selection strategy for the monotone subsequence problem. These strategies lead one in turn to a more
detailed understanding of the number of values that one actually selects.

First, one notes that it is easy to show (Samuels and Steele (1981)) that
there is a unique $y \in [0,1]$ that solves the ``equation of indifference":
$$
\widetilde{v}_{n-1}(x)= 1+ \widetilde{v}_{n-1}(y).
$$
If we denote this solution by $\alpha_n(x)$, we can use its values to determine the rule for making the sequential selections.
At the moment just before
$X_i$ is presented, we face the problem of selecting a monotone sequence from among the $n-i+1$ values $X_i, X_{i+1}, \ldots, X_n$,
and if we let $\widetilde{S}_{i-1}$ denote the last of the values $X_1, X_2, \ldots, X_{i-1}$ that has been selected so far,
then we can only select
$X_i$ if it is not greater than the most recently selected value $\widetilde{S}_{i-1}$.
In fact, one would choose to select $X_i$ if and only if it falls in the interval
$[S_{i-1}, S_{i-1}-\alpha_{n-i+1}(\widetilde{S}_{i-1})]$.
Thus, the actual
number of values selected out of the original $n$ is the random variable given by
\begin{equation*}\widetilde V_n  \stackrel{\text{def}}{=} \sum_{i=1}^n \one\, \{X_i \in [\widetilde{S}_{i-1},\, \widetilde{S}_{i-1}-\alpha_{n-i+1}(\widetilde{S}_{i-1})]\}.
\end{equation*}
By the same logic, one finds that in the sequential knapsack problem the number of values
that are selected by the optimal selection rule can be written as
\begin{equation*}\label{eq:knapsack achieved}
{V}_n \stackrel{\text{def}}{=} \sum_{i=1}^n \one\,\{X_i \in [0, \alpha_{n-i+1}(S_{i-1})]\}
\end{equation*}
where now $S_{i-1}$ denotes the capacity that remains after all of the knapsack selections have been made
from the set of values  $X_1, X_2, \ldots, X_{i-1}$ that have already been observed.
By this parallel construction we have
$$\E[{V}_n]= v_n(1)=\widetilde{v}_{n}(1)=\E[\widetilde{V}_n].$$
Moreover, one has $S_0=1= \widetilde{S}_0$,
and then see the equality of the joint distributions of the vectors
$$
(S_0, S_1, \ldots, S_{n-1}) \quad \text{and} \quad (\widetilde{S}_0, \widetilde{S}_1, \ldots, \widetilde{S}_{n-1}),
$$
since the two processes $\{S_i: 0\leq i \leq n\}$ and $\{\widetilde{S}_i: 0\leq i \leq n\}$  are (temporally non-homonomous) Markov chains
that have the same transition kernel at each time epoch.

\medskip

Theorem 1 tells us now that $\E[{V}_n ] \leq \sqrt{2n}.$ By the distributional identity of $V_n$ and $\tilde V_n$ we find indirectly that
\begin{equation}\label{eq:BRagain}
\E[\widetilde{V}_n ]  \leq \sqrt{2n} \quad \text{for all } n \geq 1.
\end{equation}
It turns out that \eqref{eq:BRagain} can be proved by a remarkable variety of methods. In particular,
Gnedin(1999) gave a direct proof where one can even accommodate a random sample size $N$
and where the upper bound is now replaced with the natural proxy $(2 \E[N])^{1/2}$.
More recently, Arlotto et al. (2015a) gave two further proofs of \eqref{eq:BRagain} as consequences of
bounds that were developed for the \emph{quickest selection problem}, a sequential decision problem that
provides a kind of combinatorial dual
to the classical sequential selection problem.

The distributional identity can also be used to make some notable inferences about the knapsack problem
from what has been discovered in the theory of sequential monotone selections. For example, by building on
the work of
Bruss and Delbaen ((2001),(2004)), 
Arlotto et al. (2015b) found that 
\begin{equation} \label{eq: Var and CLT for SMSP}
\Var[\widetilde{V}_n] \sim \frac{1}{3} \sqrt{2n} \quad \text{and } \quad
\frac
{\widetilde{V}_n -\sqrt{2n}}
{3^{-1/2} (2n)^{1/4}}
\Rightarrow {\rm N}(0,1),
\end{equation}
where $\Rightarrow$ stands for convergence in law. Thus, as a consequence of the distributional identity$~V_n\sim \widetilde V_n$ one
has the same results for the knapsack variable $V_n$ for i.i.d. $U([0,1])$ random variables. 

It is also interesting to recall here the \emph{non-sequential} (or clairvoyant) selection problem where
one studies the random variable $$L_n=\max\{k: X_{i_1}< X_{i_2}< \ldots < X_{i_k}, \, \, 1\leq i_1< i_2< \cdots <i_k \leq n \}.$$
This classic problem has a long history, beautifully told in both Romik (2014) and in Aldous and Diaconis (1999). Here the most relevant part of that story is that
Baik et al. (1999) found the asymptotic distribution of $L_n$, and, in particular, they found that one has
the asymptotic relation
\begin{equation}\label{eq:BDK99}
\E[L_n]= 2 \sqrt{n} -\alpha n^{1/6} + o(n^{1/6}) \quad \text{where } \alpha=1.77108...
\end{equation}
For the dominant part $2\sqrt{n}$ we see, as already observed in Samuels and Steele (1981), the advantage of clairvoyance is essentially  the factor $\sqrt{2}.$

\subsection{Random Tilings with different shapes}

We have given so far already several applications of Theorem 1. Now we show that the generalisation provided in  Theorem 2 has much to offer for other applications. For example, look at the following "multi-type" online problem:

\medskip
We are given a connected d-dimensional subset $S$ of $\R^d$ with Lebesgue measure $s$, and a sequence of $d$-dimensional random shapes which we would like to fit into $S$ without overlappings. If the sequence contains $n_k$ shapes of type $k$, and there are $\sigma$ types of shapes, what is the maximum number of non-overlapping shapes we can hope to fit sequentially (online) into $S$?

First, neglecting the online feature, we have a convenient upper by the BRS-equality. We define $n=n_1+n_2+\cdots+n_{\sigma},$ and let
$V_1, V_2, \cdots,V_n$ be the volumes of the $n$ random shapes. If the randomization algorithm is known for  the $\sigma$ types of shapes  and the same {\it within} each class of shapes then we can compute the distribution functions $F_k, k=1,2,\cdots,\sigma$ of their volumes which are supposed to be absolute continuous random variables $X_k$. Hence from
(5) and (6)
\begin{align} \E(N(n,s)\le \sum_{k=1}^\sigma n_k F_k(t(n,s)),~ n=n_1+n_2+\cdots+ n_\sigma, \end{align} where $t(n,s)$ solves the implicit equation \begin{align}\sum_{k=1}^\sigma n_k \int_0^{t(n,s)}\,v\, dF_k(v) =s.\end{align}
Now, coming back to online-selection, how well could we do with a skilfull online selection if we select, independently of the type of shapes all those objects which have a random volume
$V\le t(n,s)$, and fit these best possible into $S$?

\noindent As an example we study a specific tiling problem involving two types of shapes.

\subsubsection{Rectangles and ellipses}
Suppose $S$ is the unit square $S=[0,1]\times[0,1]$ and that we have two types of shapes, namely rectangles (type 1) and ellipses (type 2). Let $X, Y$ be independent $U[0,1]$ random variables, and $A, B$ be independent $U[0,1/2]$ random variables. Here $X$ and $Y$ denote the random length of the sides of a rectangle, and $A$ and $B$ the random height and width of an ellipse.  Further let $X_j, Y_j,
~j=1, 2, \cdots, n$ respectively  $A_j, B_j,
~j=1, 2, \cdots , n$ be independent versions of these.\medskip
\begin{figure}[ht]
\centering
\includegraphics[width=0.75\textwidth, angle=0]{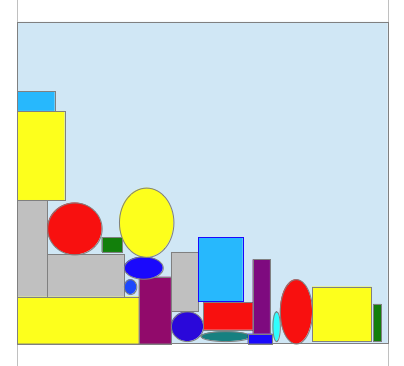}

\centerline {Figure 3}
\begin{quote}Figure 3
 displays online tiling  with  rectangles and ellipses. The square has surface $1,$ and $n_1=300$ and $n_2=150.$ Solving the BRS-equation (using Mathematica) yields $t= 0.0326\cdots$ as upper surface threshold for the rectangles and ellipses.   \end{quote}
\end{figure}
 Then we get
for the surface area distribution function $F_1(u)$ (rectangle)  for $u\in [0,1],$\begin{align} F_1(u)=P(XY\le u)=\int_0^1P\left(Y\le \frac{u}{X}|X=x\right)dx\end{align}\begin{align*}=\int_0^u dx + \int_u^1\frac{u}{x}dx=u (1-\log(u)),~~~~~~~~~~\end{align*} and similarly for the  distribution of the surface area of a random ellipse\begin{align}
F_2(u)=P\left( AB\le \frac{u}{\pi}\right)=\int_0^{1/2}P\left(B\le \frac{u}{A\pi}|A=a\right)da\end{align}
\begin{align*}=\I\{0\le u\le \pi/4\}\frac{4u}{\pi}\left (1- \log\left(
\frac{4u}{\pi}\right)\right) + \I\{\pi/4\le u\le 1\}.\end{align*}
We have $s=1$ and $\sigma=2.$ For our example (see Figure 3) let us choose $n_1=300$ and $n_2=150$, say. Putting $t:=t(450,1)$ we must solve according to the BRS-equation (29)
\begin{align}
300 \int_0^t \,u \,dF_1(u) + 150 \int_0^t\,u\, dF_2(u)~=~1,\end{align} where $F_1(u)$, respectively  $F_2(u)$ are given by (30), respectivly (31). The numerical solution (by Mathematica)
yields $t= 0.0326\cdots.$ Hence from inequality (28) $$\E(N(450,1))\le 300 F_1(t)+150 F_2(t) = 69.325\cdots.$$ 

Instead of storing $300$ rectangles and $150$ ellipses we randomized the shape $P({\rm  type}~1)=2/3, P({\rm type}~ 2)=1/3,$ and then its parameters, packing those with surface area $\le t$ independently of their type, online and just on eyesight, in the square.

In Figure 3 the upper bound $\E(N(450,1))\le 69.325\cdots$  seems surprisingly feasible. Admittedly, we were lucky with a few very small ellipses. Remember that our upper bound is a full-information
upper bound, and not only an upper bound for an online strategy. With other types of shapes the relative loss through unusable space may be higher, and for certain non-convex shapes it is likely to be much higher.
Figure 3 looks unfinished, and this was on purpose. We wanted to stop half-way, i.e. at 225 simulations. Not on purpose was that we had to stop earlier due to a legally imposed interruption, and forgotten to take note of the exact last step number (at around 150?)

The only theory we applied here was thus computing the BRS-threshold $t(n,1)$. 
It is of great help to have just a single decision function, namely the area threshold $t(n,1).$ We ask the computer to serve us sequentially only those random shapes with area below $t(n,1).$ Then each offer means "Place it!", and that's what we did (without rotations), just as they came. With the right computed threshold $t(n,s)$ it works well and is just fun to observe.

\subsection{Borel-Cantelli type problems}

For  Borel-Cantelli type problems, the BRS-inequality may also be of interest. It sometimes allow us to isolate on a set of dependent events condensed subsets on which sub-events become independent.
In our first example we use an argument in the spirit of the key-inequality (10). 

\subsubsection{Direct neigbours}
Let  $X_1, X_2, ...., X_n$ be  i.i.d. uniform  random variables on $[0,1].$  Further let $A_n$ be the event that the $(n+1)$st draw $X_{n+1}$ becomes a direct neighbour of $X_n,$ that is, $[\min\{X_n,X_{n+1}\},\max\{X_n, X_{n+1}\}]$ contains no $X_j$ with $1\le j \le n-1.$
First, we are interested to know $P(A_n~{\rm i.o.})$ as $n \to \infty,$ where "i.o." stands for infinitely-often.
Clearly $P(A_n)$ depends on where $X_n$ has landed so that the $A_n$ are not independent. Hence
$\sum_n P(A_n)=\infty$ does not imply from the  Borel-Cantelli Lemma that  $P(A_n ~i.o.)=1.$ 

However, fix now a constant $c$ with $0<c<1$ and then think of all spacings with length $\le c/n$ united in one set $A$. This takes less space than Lebesgue measure $c$ so that the difference set $B\subset [0, 1]$ has  Lebesgue measure at least $1-c,$  and this is already the simple idea: at least the Lebesgue measure $1-c$ is available for $X_{n+1}$ {\it independently} of $n,$
and within the set having at least this measure, we know that the minimal distance between neighboured points is larger than $c/n.$ 

To finish the question, let $B_{n}:=\{X_{n} {\rm ~ and~} X_{n+1} {\rm ~become~neighbours}\}.$ Clearly, the $B_{2n}$ are independent. Also, $P(B_{2n} > (1-c)c/(2n)$ since, if $X_{2n+1} \in B$, more than one half of $c/n$ is free left or right of $X_{n+1}$ for $X_{2n+2}$. Since $\{B_{2n} ~{\rm i.o.}\}\subset \{A_n~{\rm i.o.}\},$   and  the harmonic sum diverges, we have $P(A_n~{\rm i.o.})=1.$ 

\subsubsection{Spacings}
Second, let us now look at all spacings generated sequentially by the draws $X_n.$ Fix a function $k:=k(n)$ with $1\le k\le n+1$ and look at the events ${\cal K}_n=\{X_{n+1} ~{\rm falls ~in~the}~k(n)th{~\rm smallest~spacing}\}.$ Bruss et al. (1990) found criteria on $k(n)$ for $P({\cal K}_n~i.o.)=0$, respectively $P({\cal K}_n~ {\rm i.o.})=1.$ 
In a  nice more elaborate approach, Alsmeyer (1991)
improved these results by giving the complete if-and-only-if answer.

Now, what about the event of hitting infinitely often {\it one of the} $k(n)$ {\it smallest spacings}? Does it suffice to replace Alsmeyer's sum criterion by the corresponding double sum for all $j=1,2, \cdots, k(n)$?  (Since Alsmeyer  refers in his paper to decreasing order statistics we must read his criterion by replacing his $k(n)$ by $n-k(n)+1.$)
This seems intuitive, and a similar argument as above shows that this must be true. 
However, here the BRS-inequality cannot give as much as Alsmeyer's approach, which uses the strong property of uniform spacings behaving, normalized, like independent exponential variables. 

Under weaker conditions, the BRS-inequality may however give results where other approaches cannot help. For example if the BRS-sum $S_{A_(n, s_n)}$ (see (8) to (10)) "accommodates" the sum of the $k(n)$ smallest spacings for $n$ sufficiently large, and if $\E(N(n, s_n)$ stays bounded
then we {\it always} have $P({\cal K}_n~i.o.)=0.$
Moreover, if a similar "independent set" can be found as in the first question above then we are also able to find, with $\sum_{j=1}^n s_j \longrightarrow \infty$  as $n \to \infty$, the desired converse. 

\section{Further Connections and Applications}\label{se:conclusion}

Here the whole goal has been to explain, extend, and explore the
upper bound given by the BRS-inequality. The only work we know in which Theorem 2 was applied is  Steele (2016). Since, as said before, Theorem 1 has been around for quite some time (Bruss and Robertson (1991)), we know many places in the literature where Theorem 1 was used explicitly, or at least referred to. In an attempt to be somewhat complete, we give here, without giving more details, a list  in alphabetic order of those references to Theorem 1 we know:

\bigskip
Arlotto, A. and Steele, J. M. (2011),

Arlotto A., Mossel E. and Steele J. M. (2015a), (2015b)

Arlotto A.,  Wei Y., and Xie (2018)  

Boshuizen F. A. and Kertz R. P. (1995), (1999)

Bruss F. T. (2014)

Caballero J. A., Lepora N. F. and Gumey K.N. (2015)

Cain N. and Shea-Brown (2013)

Gnedin A. V. (1999)

Gnedin A. V. and Seksenbayev A. (2019)

Gribonval R., Cevher V. and Davies M. E. (2012)

Jacquemain A. (2017)

Kaluszka M., Okolevski A. and Szymanska K. (2005)

Kiwi M. and Soto A. J (2015)

Lasner Z. and  DeMille D. (2018)

Li L., Zhou H.,  Xiong S. X.,  Yang J.and Y. Mao (2019) 

Peng P. and Steele J.M. (2016)

Pohl O. (2016)

Rhee and Talagrand (1991)

Stanke M.  (2004)

Steele J.M. (2016)

Truong L V, Tan V. Y. F. (2018)

Tunc  S., Donmez M.A. and Kozat S.S. (2013)

Zhang J.J.J. (2019)

\bigskip

Theorem 1 also plays a central role in the theory
of so-called {\it resource dependent branching processes}, or RDBPs, which we will treat as our last example of applications.  The work of Bruss and Robertson (1991), instigated by the article of Coffman et al. (1987), was motivated in particular by questions arising in the context of RDBPs. Hence this is no coincidence. However, as we shall see now, it is Theorem 2 which is central for making a true progress in this theory. 

\subsection{Resource Dependent Branching processes}

A resource dependent branching process (RDPB) is a branching process trying to model the development of human populations in a way which is intended to be
as realistic as possible. 
In these RDBP's particles (individuals) have to work in order to survive and to reproduce. They receive resources from a global resource space left by their ancestors,  create new resources and claim resources for their own consumption. All these features are modelled by random variables.

 In particular, a RDBP is moreover endowed with a so-called society structure, and it is the society
who decides how resources are distributed among the individuals. It is this very feature which makes a RDBP very different from other branching process models, and this is why the methods to study them are typically also rather different. In a RDBP individuals are submitted to changing environments, but in contrast to what one usually would call a process with random random environment, the changing environment is now (primarily) the result of policies exercised by the society rather than the result of changing parameters, as for instance changing reproduction rates for the individuals.  For branching processes with random environments in the more classical sense see in particular the remarkable generality of models enabled by the unifying approach of Kersting (2019). 

\smallskip
Returning to the notion of claims used in RDBPs, it is the interplay of claims (consumption of resources) and productivity (creation of resources) which will, together with the society's policy to distribute resources, be decisive for the possible survival of the process. Namely, those individuals, who do not receive their minimum claims, are supposed  to refuse to replicate within the population. This is thus  for individuals a tool of defense, i.e. a means to have a say in policy. Finally, the  main objective of the society is to survive (first priority). The second one is to do so with a comfortable standard of living (average size of accepted claims) for the individuals. 

\smallskip
Bruss and Duerinckx (2015) proved that no society can possibly survive unless the so-called {\it weakest-first} society (w.f.) can survive forever with a strictly positive probability, where the rule of the w.f.-society is to serve the smallest  claims first. This result does not sound exciting because saving resources seems always conservative, and thus, by all means, helpful for survival. However, the logical converse (wasting resources by greedy individuals is bad for survival) should then also not be exciting, but, for RDBPs, this is a non-trivial result. Both together form the so-called Theorem of Envelopment. (See Bruss and Duerinckx (2015) for details.)

It is with the w.f.-society where the definition (1) plays a central role. Theorem 2 allows for dealing with a population which splits into several sub-populations with their own parameters. As indicated in the more elementary examples seen before, some dependencies may be fully compatible with the essence of conclusions. The final goal is to determine under which conditions sub-populations
may arrive at a certain equilibrium. 

For a long-term equilibrium between sub-populations to exist, it is a necessary condition for them to survive, and we now understand these conditions. Recall that in political decision-making (often enough online!), necessary conditions often play a more important role than sufficient conditions. We just mention that, depending on the degree of independence assumptions for the behaviour of individuals within the same sub-population, these necessary conditions
may turn as well into sufficient conditions. 

\smallskip
In particular, immigration of one or more sub-populations and its effect on equilibria can now, due to Theorem 2, be treated in one single framework  (Bruss (2018)). This is what we meant by saying that Theorem 2 will be central for any further development
of the understanding of RDBPs.
\section{Acknowledgement}
The author thanks very cordially J. Michael Steele, University of Pennsylvania, Wharton School,
for his stimulating article of 2016,  and also for so many interesting discussions we had together since then.

\section {References}

~~~~Aldous, D. and Diaconis, P. (1999), {\it~  Longest increasing subsequences: from patience
sorting to the Baik-Deift-Johansson theorem}, Bull. Amer. Math. Soc. (N.S.), 36,
413-432.

\medskip
Alsmeyer G.  (1991), {~\it Complete answer to an interval splitting problem,}
Statistics \& Probability Letters, Vol. 12,  285-287.

\medskip 
Arlotto, A. and Steele, J. M. (2011), { \it Optimal sequential selection of a unimodal
subsequence of a random sequence,} Combin. Prob. Comput., 20, 799-814.

\medskip
Arlotto A., Mossel E. and Steele J. M. (2015a), {~\it Quickest online selection of an increasing
subsequence of specified size},  Random Structures and Algorithms (eprint
arXiv:1412.7985) .
\medskip

Arlotto A., Nguyen V. V. and Steele J. M. (2015b), {~\it Optimal online selection of a monotone
subsequence: a central limit theorem}, Stochastic Process. \&Appl., Vol. 125 (9), 3596–3622.
\medskip

Arlotto A. , Wei Y.  and  Xie X. (2018) {\it An adaptive O(log n)-optimal policy for the online selection of a monotone subsequence from a random sample}, Random Structures \& Algorithms, 2018 - Wiley Online Library.

\medskip

Baik J., Deift P. and Johansson K. (1999), {~\it On the distribution of the length of the longest
increasing subsequence of random permutations}, J. Amer. Math. Soc. Vol. 12 (4), 1119–1178.

\medskip
Boshuizen F. A. and Kertz R. P. (1995), ~{\it Largest-fit selection of random sizes under a sum
constraint: comparisons by weak convergence},  in Festschrift in honour of V.S. Korolyuk (A.V. Skorokkod and Y.V. Borovskikh, eds.), 55-78.

\medskip
Boshuizen F. A. and Kertz R. P. (1999), ~{\it Smallest-fit selection of random sizes under a sum
constraint: weak convergence and moment comparisons},  Adv. in Appl. Probab. Vol. 31 (1), 178–198.

\medskip
 Bruss F. T. (2014), ~{\it Grenzen einer jeden Gesellschaft}, Jahresbericht Deutsch. Math.-Verein. Vol. 116 (3), 137-152.

\medskip 
Bruss F.T. (2018) {~\it Equilibrium Equations for Human Populations with Immigration}, arXiv:1805.01395.
 
\medskip
Bruss F.T.,  Jammalamadaka  S.R. and  Zhou X.(1990), {~\it On an
interval splitting problem}, Statistics \& Probability Letters, Vol. 10, 321-324.

\medskip
Bruss F. T. and Robertson J. B. (1991) ~{\it 'Wald's Lemma' for Sums of Order Statistics of i.i.d. Random Variables},
Adv.  Appl. Probab., Vol. 23, 612-623.

\medskip
Bruss F. T. and Delbaen F. (2001)~ {\it Optimal rules for the sequential selection of monotone subsequences of maximal length}, Stoch. Proc. \& Applic., Vol. 96, 313-342.

\medskip
 Bruss F. T. and Delbaen F. (2004), ~{\it A central limit theorem for the optimal selection process for
monotone subsequences of maximum expected length}, Stochastic Proc. \& Applic.  Vol. 114 (2), 287-
311.

\medskip
Bruss F. T. and Duerinckx M. (2015), ~{\it Resource dependent branching processes and the envelope of societies},
 Ann. of Appl. Probab., Vol. 25 (1), 324-372.
    
\medskip
Caballero J. A., Lepora N. F. and Gumey K.N. (2015),
{~\it Probabilistic Decision Making with Spikes: From ISI distributions to behavior via information gain},
~PLoS One, Vol. 10 (4); avail. at  https://doi.org/10.1371, journal.pone.0124787.
    
\medskip
Cain N. and Shea-Brown (2013),~{~\it Impact of Correlated Neural Activity on Decision-Making Performance}, Neural Computation, Vol. 25 (2), 289-327.

\medskip
Gnedin A. V. (1999), {~\it Sequential selection of an increasing subsequence from a sample of random
size}, J. Appl. Probab. Vol. 36 (4), 1074-1085.

\medskip
Gnedin A. V. and Seksenbayev A. (2019), {~\it Asymptotics and Renewal Approximation in the Online Selection of Increasing Subsequence}, arxiv 1904.11213v2.

\medskip
Coffman Jr. E. G., Flatto  L. and Weber R. R. (1987), {~\it Optimal selection of stochastic intervals
under a sum constraint}, Adv. in Appl. Probab. Vol. 19 (2), 454-473.

\medskip
Gribonval R., Cevher V. and Davies M. E. (2012), {~\it Compressible distributions for highdimensional
statistics}, IEEE Trans. Inform. Theory, Vol.  58 (8), 5016-5034.

\medskip 
Jacquemain A. (2017) {\it ~Lorenz curves interpretations of the Bruss-Duerinckx theorem for resource dependent branching processes}, arXiv:1708.01085.

\medskip
Kaluszka M., Okolevski A. and Szymanska K. (2005),
{~\it Sharp bounds for L-statistics from dependent samples of random length.} Journ. Stat. Planning
and Inference, Vol. 127, Issue 1-2, 71-89.

\medskip
Kersting G. (2019) {\it A unifying approach to branching processes in varying environments},  arXiv:1703.01960v7.

\medskip
Kiwi M. and Soto A. J (2015). {~\it Longest Increasing Subsequences of Randomly
Chosen Multi-Row Arrays.} Combinatorics, Probability and Computing, Vol. 24,  pp 254-293.

\medskip
Lasner Z. and  DeMille D. (2018)
{\it Statistical sensitivity of phase measurements via laser-induced fluorescence with optical cycling detection,} Phys. Rev. A 98, 053823;  Published online 14 Nov. 2018

\medskip Li L., Zhou H.,  Xiong S. X.,  Yang J.and Y. Mao (2019) {~\it Compound Model of Task Arrivals and Load-Aware Offloading for Vehicular Mobile Edge Computing Networks,} in IEEE Access, vol. 7, pp. 26631-26640.

\medskip
Maghsudi S. and  Hossain E. (2017). {~\it Distributed use association in energy harvesting small networks:
An exchange economy with uncertainty}, IEEE Xplore, Vol. 1 (3), 294-308.

\medskip
Peng P. and Steele J.M. (2016) {~\it Sequential selection of a monotone subsequence from a random permutation}, Proceed. Amer. Math. Society, Vol. 144, 4973-4982.

\medskip
Pohl O. (2016) {~\it Chemotaxis of sephoretic active particles and bacteria}, Doctorate thesis, Technische Universit\"at Berlin, pp 1 -126.

\medskip
Rhee W. and Talagrand  M. (1991), {~\it A note on the selection of random variables under a sum
constraint}, J. Appl. Probab. Vol. 28 (4), 919-923.

\medskip
Romik D. (2014), {~\it The Surprising Mathematics of Longest Increasing Subsequences}, Cambridge
University Press, Cambridge.

\medskip
Samuels S. M. and Steele J. M. (1981), {~\it Optimal sequential selection of a monotone sequence
from a random sample}, Ann. Probab.  Vol. 9 (6), 937-947.

\medskip
Stanke M.  (2004), {~\it Sequential selection of random vectors und a sum constraint}, J. Appl. Prob., Vol. 41 (1), 131-146.

\medskip
Steele J.M. (2016), {~\it The Bruss-Robertson Inequality: Elaborations, Extensions, and Applications}, Math. Applicanda, Vol. 44,  No 1, 3-16.

\medskip

Truong L.V.  and Tan V. Y. F. (2018)
{~\it On Gaussian MACs with Variable-Length Feedback
and Non-Vanishing Error Probabilities}, arxiv 1609.000594v3

\medskip
Tunc  S., Donmez M.A. and Kozat S.S. (2013)
{~\it Optimal investment under transaction cost: A threshold-rebalanced portfolio approach}, ~IEEE Transactions on Signal Processing, Vol. 61 (12), 3129-3142.

\medskip
Zhang J.J.J. (2019) {~\it Online Resource Allocation with Stochastic Resource
Consumption}, Department of Technology,  Stern School of Business, New York Univ., DOI: 10.13140/RG.2.2.27542.09287.

\bigskip
\noindent Author's address

\smallskip
 \noindent F. Thomas Bruss\\
Universit\'e Libre de Bruxelles\\ D\'epartement de Math\'ematique\\Campus Plaine, CP 210\\
B-1050 Brussels
\noindent (tbruss@ulb.ac.be)

\end{document}